\documentclass[12pt,leqno]{amsart}
\usepackage{amsmath}
\usepackage{amssymb}
\usepackage{hyperref}

\newtheorem{lemma}{Lemma}[section]

\newtheorem{theorem}[lemma]{Theorem}
\newtheorem{proposition}[lemma]{Proposition}
\newtheorem{definition}[lemma]{Definition}
\newtheorem{corollary}[lemma]{Corollary}
\newtheorem{example}[lemma]{Example}
\newtheorem{exercise}[lemma]{Exercise}
\newtheorem{remark}[lemma]{Remark}
\newtheorem{fig}[lemma]{Figure}
\newtheorem{tab}[lemma]{Table}

\newcommand{\bth}{\begin{theorem}}
\newcommand{\ethe}{\end{theorem}}
\newcommand{\bre}{\begin{remark}\em }
\newcommand{\ere}{\end{remark}}
\newcommand{\ble}{\begin{lemma}}
\newcommand{\ele}{\end{lemma}}

\newcommand{\bde}{\begin{definition}}
\newcommand{\ede}{\end{definition}}

\newcommand{\bco}{\begin{corollary}}
\newcommand{\eco}{\end{corollary}}

\newcommand{\bpr}{\begin{proposition}}
\newcommand{\epr}{\end{proposition}}

\newcommand{\bexer}{\begin{exercise}}
\newcommand{\eexer}{\end{exercise}}

\newcommand{\bexam}{\begin{example}}
\newcommand{\eexam}{\end{example}}

\newcommand{\bfi}{\begin{fig}}
\newcommand{\efi}{\end{fig}}

\newcommand{\btab}{\begin{tab}}
\newcommand{\etab}{\end{tab}}

\newcommand{\beao}{\begin{eqnarray*}}
\newcommand{\eeao}{\end{eqnarray*}\noindent}

\newcommand{\beam}{\begin{eqnarray}}
\newcommand{\eeam}{\end{eqnarray}\noindent}

\newcommand{\beqq}{\begin{equation}}
\newcommand{\eeqq}{\end{equation}\noindent}

\newcommand{\bce}{\begin{center}}
\newcommand{\ece}{\end{center}}

\newcommand{\barr}{\begin{array}}
\newcommand{\earr}{\end{array}}

\newcommand{\bdis}{\begin{displaymath}}
\newcommand{\edis}{\end{displaymath}\noindent}

\begin{document}

\title{Characterization of tails through hazard rate and convolution closure properties }

\author[A. G. Bardoutsos]{A. G. Bardoutsos}
\thanks {Dept. of Statistics and Actuarial - Financial Mathematics, 
University of the Aegean\\ Karlovassi, GR-83 200 Samos, Greece\\ e-mail: sasm10007@sas.aegean.gr}

\author[D. G. Konstantinides]{D. G. Konstantinides} 
\thanks{Dept. of Statistics and Actuarial - Financial Mathematics, 
University of the Aegean\\ Karlovassi, GR-83 200 Samos, Greece\\ Tel.  +3022730-82333, Fax. +3022730-82309, e-mail: konstant@aegean.gr}

\begin{abstract}
We use the properties of the Matuszewska indices 
to show asymptotic inequalities for hazard rates. 
We discuss the relation between membership in the classes of 
dominatedly or extended rapidly varying tail distributions and 
corresponding hazard rate conditions. 
Convolution closure is established for the class of 
distributions with extended rapidly varying tails.
\end{abstract}


\maketitle

\section{Introduction}
\setcounter{equation}{0}

In this paper we intend to discuss Pitman's criterion for 
subexponentiality (see \cite[Theorem 2]{pitman:1980}).
Some extensions of previous results about the characterization of distribution classes through their hazard rates
appeared as  byproducts. 
The motivation was the need for understanding and calculating the
monotonicity condition required in these theorems. 
The ultimate goal is to substitute the monotonicity property with some limit relation.

Consider the Lebesgue convolution for densities $f_1$ and $f_2$ on $[0\,,\infty)$:
\beao
f_1 \star f_2(x)=\int_{0}^{x}f_1(y)f_2(x-y)dy\,,
\eeao
and the convolution formula for the corresponding distributions
\beao
\overline{F_1 \ast F_2}(x)=\overline{F}_2(x)+\int_{0}^{x}\overline{F}_1(x-y)dF_{2}(y)\,,
\eeao
where $\overline{F}(u)=1-F(u)$ denotes the right tail of any distribution $F$.

For $u>1\,,$ write $\overline{F}_{\star}(u):=\liminf_{x \rightarrow
  \infty}\overline{F}(ux)/\overline{F}(x)$ and
$\overline{F}^{\star}(u):=\limsup_{x \rightarrow
  \infty}\overline{F}(ux)/\overline{F}(x)$. We write $m(x)
\sim g(x)$ as $x \rightarrow \infty$ for the limit relation $\lim_{x
  \rightarrow \infty} m(x)/g(x) =1$ and introduce 
the following classes of distributions $F$:

\begin{enumerate}

\item We say that $F$ has extended rapidly varying tail\,, if $\overline{F}^{\star}(u)<1$\,, for some $u>1$. We write $F \in \mathcal{E}.$

\item  $F$ has a subexponential distribution, if $ {\overline{F^{2\ast}}(x)}\sim 2{\overline{F}(x)}$. We write $F \in \mathcal{S}$.
\item $F$ has a long tail if $\overline{F}(x-y)\sim \overline{F}(x)$ for $y \in (-\infty\,, \infty)$. We write $F \in \mathcal{L}$. 
\item $F$ has dominatedly varying tail, if $\overline{F}_{\star}(u)>0$ for all (or eq. for some) $u>1$ or, equivalently, $\overline{F}^{\star}(u)<\infty$, for all (or eq. for some) $0<u< 1$. We write $F \in \mathcal{D}$.
\end{enumerate}
Recall that for a positve function $g$ on
$(0,\infty)$ the upper Matuszewska index $\gamma_g$ is defined as the
infimum of those values $\alpha$ for which there exists a constant $C$
such that for each $U>1$, as $x\to\infty$,
$$
g(ux)/g(x)\le C(1+o(1))\,u^\alpha\quad\mbox{uniformly in $u\in
[1,U]$,}
$$
and the lower  Matuszewska index $\delta_g$ is defined as the
supremum of those values $\beta$ for which, for some $D>0$ and all $U>1$, as $x\to\infty$,
$$
g(ux)/g(x)\ge D(1+o(1))\,u^\beta\quad\mbox{uniformly in $u\in
[1,U]$.}
$$
The classes $\mathcal{D}$ and $\mathcal{E}$ are linked  to the
Matuszewska indices of the tails $\overline F$ (see
\cite{Cline-Sam:1994}).
For any distribution $F$ on $(0,\infty)$ with infinite support, $F \in 
\mathcal{D}$ if and only if $\gamma_{\overline{F}}<\infty$, and $F \in
\mathcal{E}$ if and only if $\delta_{\overline{F}}>0$. In what
follows, we always assume that $F$ has a  positive Lebesgue density
$f$. Then the following relations hold (see \cite[Theorem 2.1.5]{bingham:goldie:teugels:1987}):
$$
\gamma_{f}=\inf\left\{-\log{f_{\star}(u)}/\log{u} :u>1\right\}=-\lim_{u\rightarrow \infty}\log{f_{\star}(u)}/\log{u}\,,
$$
where $f_{\star}(u)=\liminf_{x \rightarrow \infty} f(ux)/f(x)$\,, and 
$$
\delta_{f}=\sup\left\{-\log{f^{\star}(u)}/\log{u} :u>1\right\}=-\lim_{u\rightarrow \infty}\log{f^{\star}(u)}/\log{u}\,,
$$
where $f^{\star}(u)=\limsup_{x \rightarrow \infty} f(ux)/f(x)$. 
Using the Matuszewska indices, one can establish 
Potter-type inequalities for $f$; see \cite[Proposition
2.2.1]{bingham:goldie:teugels:1987}. For example,
if $\gamma_f<\infty$, then for every $\gamma > \gamma_f$ there exist constants $C^{\prime}(\gamma)\,,\;x^{\prime}_0=x^{\prime}_0(\gamma)$ such that  
\begin{equation} \label{M 3.1}
f(y)/f(x)\ge C^{\prime}(\gamma) (y/x)^{-\gamma}\,, \quad y\ge x \ge x_0^{\prime}.
\end{equation}
If $\delta_f>-\infty$ then for every $\delta < \delta_f$ there exist constants $C(\delta)\,,\;x_0=x_0(\delta)$ such that 
\begin{equation} \label{M 3.2}
f(y)/f(x)\le C(\delta) (y/x)^{-\delta} \,, \quad y\ge x \ge x_0.
\end{equation}
In what follows, we say that
the distributions $F_1$ and $F_2$ are  max-sum equivalent, if 
\beao
\lim_{x \rightarrow \infty}\overline{F_1 \ast F_2}(x)/(\overline{F}_1(x)+\overline{F}_2(x))=1.
\eeao

\section{The class of distributions with extended rapidly varying
  tails}
We say that the distribution $F$ on $(0,\infty)$ is heavy tailed if
$\int_{0}^{\infty}e^{sx}dF(x)=\infty$ for all $s>0$ and light tailed, otherwise. 
The class $\mathcal{E}$ contains both light and heavy tailed
distributions. For example, the exponential and Pareto distributions are
members of $\mathcal{E}$. Moreover, the class $\mathcal{E}$ is not
closed under max-sum equivalence. For example, for exponential $F_i=F$
exponential with parameter $\lambda$ we have 
$
\lim_{x \rightarrow \infty}\overline{F \ast F}(x)/(2\overline{F}(x))=\infty.
$

In what follows, we will need the hazard rate 
$h(x)=f(x)/\overline F(x)$ for any distribution $F$ on
$(0,\infty)$ with positive density $f$. 
We also write
\beao
M_1= \liminf_{x \rightarrow \infty} xh(x)\quad\mbox{and}\quad   M_2= \limsup_{x \rightarrow \infty} xh(x)\,.
\eeao
Whenever we consider a sequence $F_i$, $i=1,2,\ldots,$ of such distributions, we will
use the corresponding symbols $h_i,M_1^i,M_2^i$.

We say that a density has bounded increase 
if $\delta_f>-\infty$; see 
\cite[Definition, page 71]{bingham:goldie:teugels:1987}. 
Most of the densities of interest in statistics and probability theory
satisfy this condition, e.g. the gamma and Weibull densities. 

Under the assumption of an eventually 
non-increasing density (such that $f(y)\le f(x)$, for all $x\ge y\ge x_0$), the following equivalences were established:
\begin{enumerate}
\item $F \in \mathcal{E}$ if and only if $M_1>0$; see \cite[Proposition 6]{konstantinides:2008}.
\item $F \in \mathcal{D}$ if and only if $F \in \mathcal{D}\cap \mathcal{L}$ if and only if $M_2<\infty$; see \cite[Corollary 3.4]{kluppelberg:1988}. 
\end{enumerate}

We generalize these results by substituting 
the condition of an eventually non-increasing density $f$
by the assumption that $f$ has bounded increase. This allows one to avoid 
the verification of the monotonicity property of $f$, and restricts
one to the calculation of  $\delta_f$ through $f^{\star}(u)$. 
\bth \label{pr 1.7}
Assume $F$ is a distribution supported on $(0,\,\infty)$ 
with positive Lebesgue density $f$ 
such that $f$ has bounded increase. Then
$F \in \mathcal{E}$ if and only if  $M_1>0$.
\ethe
\begin{proof}
By assumption, we have $\delta_f>-\infty$ and therefore inequality
\eqref{M 3.2} holds for $\delta<\delta_f$. Let us start with the 
converse assertion, i.e. we assume $M_1>0$.  We write for any $u>1$
\beam \label{G_1}
\overline{F}(ux)/\overline{F}(x)=\overline{F}(ux)/\left(\int_{x}^{ux}f(t)dt+\overline{F}(ux)\right)\,.
\eeam
Now the relation \eqref{M 3.2} implies that
\beao
f(ux)/f(t)\le C(\delta)(ux/t)^{-\delta} \,, \quad ux \ge t >x\ge x_{0}\,.
\eeao
Then integration yields 
\beam \label{y}
\int_{x}^{ux} f(t)dt\ge (f(ux)/C(\delta))\,(ux)^{\delta}\int_{x}^{ux}t^{-\delta}dt :=K(\delta,\,u)x\,f(ux)\,.
\eeam 
We substitute this lower bound into \eqref{G_1} to obtain
\beao
\dfrac{\overline{F}(ux)}{\overline{F}(x)}\le\dfrac{\overline{F}(ux)}{K(\delta,u)x f(ux)+\overline{F}(ux)}=\dfrac{1}{K(\delta,u)x h(ux)+1}\,.
\eeao
Since $M_1>0$ the latter relation  implies $F \in \mathcal{E}$. 

Now we show the direct implication. From \eqref{G_1} we obtain
\beam \label{b_1}
\overline{F}(ux)/\overline{F}(x)=1-\int_{x}^{ux}\left(f(t)/\overline{F}(x)\right)\,dt\,,
\eeam
and from \eqref{M 3.2} we have
\beao
f(t)/f(x)\le C(\delta)(t/x)^{-\delta}, \quad t\ge x\ge x_0\,. 
\eeao
An approach similar to the one in the inequality \eqref{y} yields for
some constant $K'(\delta,u)$
\beam \label{b_3}
\int_{x}^{ux} f(t)dt \le C(\delta)\, f(x) x^{\delta}\int_{x}^{ux}t^{-\delta}dt \le  K'(\delta,u)\,xf(x)\,.
\eeam
Hence, from (\ref{b_1}) and (\ref{b_3}) it follows that
\beao
\overline{F}(ux)/\overline{F}(x)\ge 1-K'(\delta,u) xh(x)\,.
\eeao
Since $F\in {\mathcal E}$ there exists $u>1$ such that
\beao
M_1=\liminf_{x\rightarrow \infty} xh(x)\ge (K'(\delta,u))^{-1}\liminf_{x \rightarrow \infty}(1-\overline{F}(ux)/\overline{F}(x))>0 \,.
\eeao
This finishes the proof.
\end{proof}

We say that a positive Lebesgue density $f$ is extended rapidly
varying if $\delta_f>1$. In the following result we prove that this
property and bounded increase of $f$ imply $F \in \mathcal{E}$.  We find an asymptotic lower bound for the hazard rate using the lower Matuszewska index. 

\bpr \label{prop 2.2}
If $f$ has bounded increase with $\delta_f>1$ then $F \in \mathcal{E}$
and for any $\delta \in(1,\delta_f)$ there are positive constants
$x_0\,, \; C(\delta)$, 
defined in  \eqref{M 3.2}, such that for all $x\ge x_0$:
\beam \label{eq.2.6}
xh(x)\ge (\delta-1)/C(\delta)\,.
\eeam
\epr
\begin{proof} 
The inequality \eqref{M 3.2} is implied by the assumption $\delta_f>1$
for $\delta<\delta_f$. Further, we integrate \eqref{M 3.2} for
$\delta\in (1,\delta_f)$:
\beam \label{eq 3.7}
(h(x))^{-1}=\int_{x}^{\infty}f(y)/f(x) dy \le C(\delta) x^{\delta} \int_{x}^{\infty} y^{-\delta} dy=C(\delta)x/(\delta-1)\,.
\eeam
This proves inequality \eqref{eq.2.6}. The latter relation immediately
implies that $M_1>0$ and therefore, by Theorem \ref{pr 1.7},
$F\in {\mathcal E}$.
\end{proof}
We examine the convolution closure of distributions with extended rapidly varying distributions. 
We write $h_{F_1 \ast F_2 }$ for the hazard rate  of $F_1 \ast F_2$ and 
\beao 
M_1^{(1,2)}:=\liminf_{x \rightarrow \infty} x\,h_{F_1 \ast F_2 }(x)\,.
\eeao
We need the fact that if $F_i\in \mathcal{E}$ then for every
$0<\delta<\delta_{\overline{F}_i}$ there exist constants $x_0^i =
x_0^i(\delta)$ and $C_i(\delta)$ such that the following Potter-type
inequality holds (see \cite{Cline-Sam:1994})
\beao
\overline{F_i}(x)/\overline{F}_i(y)\le C_i(\delta)
(x/y)^{-\delta}\,, \quad x\ge y\ge x_0^i\,,\quad
i=1,2\,.
\eeao
Choosing $y=x_0^i$, the latter relation implies that there exist
constants $\Delta_i(\delta)$ such that
\beam \label{hg}
\overline{F}_i(x) \le \Delta_i(\delta) x^{-\delta}\,,\quad x\ge
x_0^i\,, \quad i=1,2\,.  
\eeam
\bth \label{th 2.5}
Assume that $F_1\,, F_2 \in \mathcal{E}$ with positive Lebesgue
densities on $(0,\infty)$ and that the following conditions hold:
\begin{enumerate} 
\item The density $f_1$ has bounded increase with $\delta_{f_1}>0\,,$
\item $\delta_{\overline F_1}<\delta_{\overline F_2}$ and 
$\liminf_{x\rightarrow
    \infty}x^{{\delta}}\,\overline{F}_1(x)> 0$ for some
  $\delta\in \left[\delta_{\overline F_1 },\delta_{\overline F_2
  }\right)$.
\end{enumerate}
Then $F_1 \ast F_2 \in \mathcal{E}$. 
\ethe
\begin{proof}
We start by proving $M_1^{(1,2)}>0$. We have
\beao
x\,h_{F_1 \ast F_2 }(x) =\dfrac{xf_1(x)}{\overline{F}_1(x)} \dfrac{\overline{F}_1(x)}{\overline{F_1 \ast F_2}(x)}\dfrac{f_1 \star f_2(x)}{f_1(x)}\,.  
\eeao 
By assumption $1.$,  inequality \eqref{M 3.2} applies for
$\delta_0<\delta_{f_1 }$. Therefore for $x\ge x_0$, 
\beao
\dfrac{f_1 \star f_2(x)}{f_1(x)}&=&\int_{0}^{x}\dfrac{f_1(y)}{f_1(x)}\,f_2(x-y)dy \ge  \int_{x_0}^{x} 
\dfrac{1}{C(\delta_0)} \left(\dfrac{y}{x}\right)^{-\delta_0} \,f_2(x-y)dy \\ [2mm]
&\ge& F_2(x-x_0)/C(\delta_0)\,.
\eeao
Now, from \eqref{hg} for every $\delta_i \in (0,\delta_{\overline F_i})$,
$i=1,2$, and sufficiently large $x$, some constant $\Delta>0$,\,
\beao
\overline{F_1 \ast F_2}(x) &\le&
\overline{F}_1(x/2)+\overline{F}_2(x/2)\\
&\le & \Delta_1(\delta_1) (x/2)^{-\delta_1}+\Delta_2(\delta_2) (x/2)^{-\delta_2} \le \Delta\left(x^{-\delta_1}+x^{-\delta_2}\right)\,.
\eeao
We conclude that there exist $\delta_2 \in \left[\delta_{\overline F_1 },\delta_{\overline F_2}\right)$ such that
\beao
M_1^{(1,2)}&\ge &
\liminf_{x \rightarrow \infty} \dfrac{1}{\Delta\, C(\delta_0)}
\dfrac{xf_1(x)}{\overline{F}_1(x)}
x^{\delta_2}\overline{F}_1(x)
\ge   \dfrac{M_1^1}{\Delta\, C(\delta_0)} \liminf_{x\to\infty} x^{\delta_{2}}\overline{F}_1(x)>0\,.  
\eeao
In the last step we used Theorem \ref{pr 1.7} for  $M_1^1>0$ and 
assumption 2. 
Another  application of Theorem \ref{pr 1.7} yields the result. 
\end{proof}
To verify that the two conditions do not contradict each other
we consider two Pareto distributions with  tails 
$\overline F_i(x)=x^{-a_i}$, $i=1,2$, $x\ge 1$. Choose $a_1=2$ and 
$a_2=3$. Then it is easy to see 
that $\delta_{\overline{F}_i}=a_i$, $i=1,2$. Therefore 
$F_i\in \mathcal{E}$ and $f_i^{\star}(u)=u^{-(a_i+1)}$. Hence
condition $1.$ holds with $\delta_{f_1}=3$, $f_1$ has bounded
increase and condition 2. is satisfied for $\delta=2$.
We conclude that the conditions of the theorem hold. 
The tail of the Pareto distribution $F_i$ belongs to the class $\mathcal{R}_{-a_i}$
of regularly varying functions with index $-a_i$, i.e. 
$\overline{F_i}(ux)\sim u^{-a_i}\overline{F_i}(x)$ for each $u>0$. 
From \cite[Lemma 1.3.1]{embrechts et.al} we know that if 
$\overline F_i \in \mathcal{R}_{-a_i}$, $i=1,2$,  then 
$\overline{F_1\ast F_2}\in \mathcal{R}_{-\min(a_1,a_2)}$. This fact
and the definition of regular variation immediately imply $F_1\ast
F_2\in {\mathcal E}$.

\section{The case of subexponential distributions}

In this section we present some results on the characterization 
of the classes $\mathcal{S}$ and $\mathcal{D}\cap\mathcal{L}$ through 
hazard rates. In the following result 
we provide an inequality for the hazard rate which is useful for
characterizing membership in the class $\mathcal{D}\cap \mathcal{L}$. 

\bpr \label{lem 2.1} Assume that $F$ has a positive Lebesgue
  density on $(0,\infty)$ and $\gamma_f<\infty$.
Then  $F \in \mathcal{D} \cap \mathcal{L}$, and for any $\gamma >
\gamma_f$ there are positive constants $x_0^{\prime}\,, \;
C^{\prime}(\gamma)$, defined in \eqref{M 3.1}, such that for all $x\ge
x_0^{\prime}$ and $\lambda>1$: 
\beam\label{aux}
x h(x)\le C^{\prime}(\gamma)\,V(\lambda, \gamma)\,,
\eeam
where
\beao
V(\lambda, \gamma)=\left\{
\begin{array}{ll}
(\lambda^{-\gamma+1}-1)/(-\gamma+1),& \mathrm{if}\; \gamma\neq 1,\\[2mm]
\log\lambda, & \mathrm{if}\; \gamma=1.
\end{array}\right.
\eeao
\epr

\begin{proof}
Since $\gamma_f<\infty$, \eqref{M 3.1} yields
\beao
\int_{x}^{\infty}f(y)/f(x)dy &\ge& \int_{x}^{\lambda x}f(y)/f(x) dy \ge C^{\prime}(\gamma)\, x^{\gamma} \int_{x}^{\lambda x} \nonumber y^{-\gamma}dy\,.
\eeao
Then \eqref{aux} holds, $M_2<\infty$ and
from \cite[Theorem 3.3]{kluppelberg:1988} we obtain $F \in \mathcal{D}\cap \mathcal{L}$. 
\end{proof}

In the next theorem we generalize the statement from \cite[Corollary 3.4]{kluppelberg:1988} by substituting 
the condition of an eventually non-increasing density $f$
by the assumption that $f$ has bounded increase. This allows one to avoid 
the verification of the monotonicity property of $f$, and restricts
one to the calculation of  $\delta_f$ through $f^{\star}(u)$.  

\bth \label{D_0}
Assume that $F$ is supported on $(0,\,\infty)$ with a positive
Lebesgue density $f$ which has  bounded increase. 
Then the following statements are equivalent:\beao
\mbox{1. $F \in \mathcal{D}$\qquad 
2. $F \in \mathcal{D}\cap \mathcal{L}$\quad and\quad   
3. $M_2< \infty$\,.}
\eeao 
\ethe
\begin{proof}
$(1)\Rightarrow (3)$.
We start by observing that
\beam \label{56}
\overline{F}(x/2)/\overline{F}(x)=1+\int_{x/2}^{x}f(y)/\overline{F}(x)\,dy\,.
\eeam 
Since $f$ has bounded increase \eqref{M 3.2} applies for $x\ge y \ge
x/2$ and sufficiently large $x$. Hence exists a constant
$\Lambda(\delta)$ such that for large $x$,
\beam\label{57}
\int_{x/2}^{x}f(y)dy\ge (x^{\delta}f(x)/C(\delta))\int_{x/2}^{x}y^{-\delta}dy=xf(x)\Lambda(\delta)\,.
\eeam
The inequalities \eqref{56} and \eqref{57} imply
\beao
xh(x)\le [\Lambda(\delta)]^{-1}(\overline{F}(x/2))/\overline{F}(x)-1)\,.
\eeao
Now, from the assumption $F \in \mathcal{D}$ we obtain $M_2< \infty$. \\
$(3) \Rightarrow (2)$ follows from \cite[Theorem
3.3]{kluppelberg:1988} and 
$(2)\Rightarrow (1)$ is trivial.
\end{proof}

In \cite[Theorem 2]{pitman:1980} necessary and sufficient conditions
for membership in $\mathcal{S}$ were presented. 
In \cite[Theorem 3.6]{kluppelberg:1988} a corresponnding result for
a the important subclass ${\mathcal S}^\ast$ of $\mathcal{S}$ was given; see also \cite[Corollary
3.8]{kluppelberg:1988}. The previously mentioned results require that the hazard rate be
eventually monotone (such that $h(y)\le h(x)$ for all $y\ge x\ge x_0$).
However, a verification of this monotonicity condition is in general
not straightforward. In the next result we prove the statement of
\cite[Theorem 2]{pitman:1980} under the assumption $\delta_h>0$
which might be checked easier. 

Recall the notion of hazard function $H(x):=-\ln \overline{F}(x)$ with
the convention $H(\infty)=\infty$. An application of \cite[Proposition
2.2.1]{bingham:goldie:teugels:1987} yields a Potter-type inequality
for $g(x)=(h(x))^{-1}$: if $\delta_h>0$ then for
$0<\delta < \delta_h$ there exist  constants $C(\delta)\,,\;x_0$ such that 
\begin{equation} \label{H 3.2}
h(y)/h(x)\le C(\delta) (y/x)^{-\delta}\,, \quad y\ge x \ge x_0 \,.
\end{equation}
If $\delta_h>0$ we say that the hazard rate $h$ has positive decrease; see \cite[Definition, page~71]{bingham:goldie:teugels:1987}.
\bth \label{pitman0}
Let $F$ be a distribution on $(0,\infty)$ with positive Lebesgue
density $f$. Assume that the hazard rate $h$ has positive decrease. Then $F\in \mathcal{S}$ if and only if 
\beam \label{Pitman}
\lim_{x\rightarrow \infty} \int_{0}^{x}\exp\{\kappa y h(x)-H(y)\}h(y)dy=1\quad\mbox{for every $\kappa>0$.}
\eeam 
\ethe
\begin{proof}  
If $F$ has a  positive Lebesgue density the hazard function $H$ is
differentiable and we obtain 
\beao
\overline{F^{2\ast}}(x)\big/ \,\overline{F}(x)-1
&=&\int_{0}^{x} \exp\left\{H(x)-H(x-y)-H(y)\right\} h(y)dy \nonumber\\
&=&\int_{0}^{x/2} \exp\left\{H(x)-H(x-y)-H(y)\right\} h(y)dy \nonumber\\
&&+\int_{0}^{x/2} \exp\left\{H(x)-H(x-y)-H(y)\right\}h(x-y)dy \nonumber \\
&=&I_1(x)+I_2(x)\,, 
\eeao
We start by showing the 
converse  implication: \eqref{Pitman}  implies that $I_1(x) \rightarrow 1$ and
$I_2(x) \rightarrow 0$, hence $F$ is subexponential.
For $y\le x/2$ there exists  $\xi\in (x-y\,,x)$ such that 
$y h(\xi)=H(x)-H(x-y)$. Then
\beam \label{W_11}
x>\xi>x-y\ge x/2\ge y.
\eeam
An application of  \eqref{H 3.2} yields for large $x$ and $\delta<\delta_h$
\beam \label{T2}
\max(h(\xi)/h(x/2),
h(x)/ h(\xi))\le C(\delta)\,.
\eeam
Hence for any $x,\,y$ satisfying \eqref{W_11}, 
\beam \label{W_45}
y\,h(x)/ C(\delta)\le H(x)-H(x-y) \le  
C(\delta)\,y \,h( x/2)\,, 
\eeam
Since $H(x)-H(x-y)\ge 0$ we have the trivial bound $I_1(x)\ge F(x/2)$ and
from the right inequality in \eqref{W_45} we conclude that
\beam\label{eq:0}
F(x/2)\le I_1(x)\le \int_{0}^{x/2} \exp\{C(\delta)\,y h(x/2)-H(y)\} h(y)\,dy\,.
\eeam
Together with
\eqref{Pitman} this implies the desired relation $I_1(x)\rightarrow 1$.

It remains to show $I_2(x) \rightarrow 0$.  
From \eqref{W_45} we obtain 
\beao
{I}_2(x)\le \Big(\int_0^{x_0}+\int_{x_0}^{x/2}\Big)\exp\{ C(\delta)\,y h(x/2)-H(y)\}h(x-y)\,dy\,.
\eeao
The first integral converges to zero as $x\to\infty$ since
$h(x)\rightarrow 0$. 
By \eqref{W_11} and \eqref{H 3.2} there exists some $x_0$ such that 
$h(x-y)/ h(y) \le C(\delta)$, $x_0 \le y \le x/2$, and therefore, up
to a constant multiple, the second integral is bounded by the
right-hand expression in \eqref{eq:0} which converges to 1. Moreover,
the integrand in the right-hand expression of \eqref{eq:0} converges
for every $y$ as $x\to\infty$. The integrand in the second integral
above converges to zero for every $y$. An application of   
Pratt's lemma \cite[Theorem 1]{pratt:1960} shows that the second
integral
converges to zero as $x\to\infty$. Hence $I_2(x)\to0$ and the converse
implication of the result is proved.

For the direct part, assuming that $F$ is subexponential, 
we obtain from \eqref{W_11} and \eqref{T2} that 
\beao
\overline{F^{2\ast}}(x)\big/ \, \overline{F}(x)-1\ge I_1(x)
\ge\int_{0}^{x/2} \exp\left\{-H(y)\right\}h(y)dy 
=\int_{0}^{x/2} f(y)dy\,. 
\eeao
Since the left-and right-hand sides converge to 1 as $x\to\infty$ the
proof of \eqref{Pitman} is complete.
\end{proof}

We introduce the class
$\mathcal{A}=\mathcal{S}\cap \mathcal{E}$ of heavy tailed
distributions on $(0,\infty)$. In what
follows, we will frequently use the relation
$\mathcal{D}\cap \mathcal{A}=\mathcal{D}\cap
\mathcal{L}\cap \mathcal{E}$, which is a consequence of \cite[Theorem 1]{goldie:1978} 
and the definition of the class $\mathcal{A}$.
Also notice that the inclusion 
$\mathcal{D}\cap \mathcal{A}\subset
\mathcal{D}\cap\mathcal{L}$ is strict since 
the distributions with slowly varying
tail are not contained in $\mathcal{E}$, but belong to
$\mathcal{D}\cap \mathcal{L}$. 
Next we show that $\mathcal{D}\cap \mathcal{A}$ is closed under convolution.
\bpr Assume
$F_i \in  \mathcal{D}\cap \mathcal{A}$, $i=1,2$. 
Then $F_1\ast F_2 \in  \mathcal{D}\cap \mathcal{A}$ and 
\beam \label{adms}
\overline{F_1\ast F_2}(x) \sim \overline{F}_1(x)
+\overline{F}_2(x)\,,\quad x\to\infty\,.
\eeam
\epr
\begin{proof}
Relation \eqref{adms} follows from 
$F_i\in \mathcal{D}\cap\mathcal{A} \subset
\mathcal{D}\cap\mathcal{L}$, $i=1,2,$ and 
\cite[Theorem 2.1]{cai:tang:2004}.  
Furthermore from \cite[Proposition 2]{embrechts:goldie:1980} (or
\cite[Theorem 2.1]{cai:tang:2004}) 
we conclude 
$F_1\ast F_2 \in \mathcal{D}\cap\mathcal{L}$. Let us use relation \eqref{adms} to show that  $F_1\ast F_2 \in  \mathcal{E}$. Then for some $u>1$ we have
\beao
\limsup_{x \rightarrow \infty}\dfrac{\overline{F_1\ast F_2}(ux) }{\overline{F_1\ast F_2}(x)}=\limsup_{x\rightarrow \infty}\dfrac{ \overline{F}_1(ux)+ \overline{F}_2(ux)}{ \overline{F}_1(x)+ \overline{F}_2(x)}\le \max\left\{\overline{F}^{\star}_1(u)\,,\;\overline{F}^{\star}_2(u)\right\}.
\eeao
Taking into account that $F_i\in \mathcal{E}\,,\;i=1,\,2$ we obtain
the result. This finishes the proof.
\end{proof}

In the next statement a characterization  of the class
$\mathcal{D}\cap \mathcal{A}$ with respect to the hazard rate and the
limits $\overline{F}_{\star}(u)$ and  $\overline{F}^{\star}(u)$ for
all $u>1$ is presented and in this way a generalization of
\cite[Theorems 3.3 and 3.7]{konstantinides:2008} is provided. The generalization is achieved by substituting 
the condition of an eventually non-increasing density $f$
by the assumption that $f$ is of bounded increase. This allows one to avoid 
the verification of the monotonicity property of $f$, and restricts
one to the calculation of  $\delta_f$ through $f^{\star}(u)$.  

\bco
Assume $F$ is a distribution supported on $(0,\,\infty)$ 
with positive Lebesgue density $f$ 
such that $f$ has bounded increase. Then $F \in \mathcal{D}\cap \mathcal{A}$ if and only if one of the following statements holds
\begin{enumerate}
\item $0<M_1\le M_2 <\infty\,,$ 
\item $0<\overline{F}_{\star}(u)\le\overline{F}^{\star}(u)<1\,.$ 
\end{enumerate}
\eco
\begin{proof}

$1.$ Let us begin with the direct implication. From the assumption $F\in \mathcal{D}$ and Theorem \ref{D_0} we obtain $M_2<\infty$. From the assumption $F \in \mathcal{E}$ and from Theorem \ref{pr 1.7} we find $M_1>0$. Next, for the inverse part we apply directly Theorems \ref{pr 1.7} and \ref{D_0}.  

$2.$ The condition  $F \in \mathcal{D}\cap\mathcal{A}$ is equivalent to  $F \in \mathcal{D}\cap\mathcal{E}$. Because of $F \in \mathcal{D}$ from Theorem \ref{D_0} we obtain $F \in \mathcal{D}\cap\mathcal{L}$. Hence $F \in \mathcal{D}\cap\mathcal{A}$. This completes the proof.
\end{proof}

\bco
Let $F$ be a distribution on $(0,\infty)$ with positive Lebesgue
density $f$. Assume that the hazard rate $h$ has positive decrease. Then $F\in \mathcal{A}$ if and only if $M_1>0$ and
\beam \label{Pitman2}
\lim_{x\rightarrow \infty} \int_{0}^{x}\exp\{\kappa y h(x)-H(y)\}h(y)dy=1\quad\mbox{for every $\kappa>0$}\,,
\eeam 

\eco
\begin{proof}
From the decreasing property of the tail $\overline{F}$ we obtain $f(ux)/f(x)\le h(ux)/h(x)$ for every $u>1$. Hence $\delta_f\ge\delta_h>0$. From the last inequality and Theorems \ref{pr 1.7} and \ref{pitman0} we conclude the result. This finishes the proof.
\end{proof}

\vskip0.5cm

\noindent \textbf{Acknowledgments.} The authors would like to thank
Thomas Mikosch 
for his comments that improved significantly the material of the paper
and S\o ren Asmussen for his advice to work on this topic.

\end{document}